\documentclass[a4paper]{amsart}
\usepackage{graphicx}
\usepackage{amsmath}
\usepackage{amssymb}
\usepackage{amsthm}
\usepackage{oldgerm}
\usepackage[all]{xy}
\usepackage[latin1]{inputenc}

\newcommand{\Hom}{\operatorname{Hom}\nolimits}
\newcommand{\End}{\operatorname{End}\nolimits}

\newcommand{\Ker}{\operatorname{Ker}\nolimits}
\newcommand{\Coker}{\operatorname{Coker}\nolimits}

\newcommand{\Ext}{\operatorname{Ext}\nolimits}
\newcommand{\op}{\operatorname{op}\nolimits}

\newcommand{\modu}{\operatorname{mod}\nolimits}

\newcommand{\add}{\operatorname{add}\nolimits}

\newtheorem{theorem}{Theorem}[section]

\newtheorem{lemma}[theorem]{Lemma}
\newtheorem{proposition}[theorem]{Proposition}
\newtheorem{corollary}[theorem]{Corollary}
\theoremstyle{definition} \newtheorem*{definition}{Definition}

\begin{document}
\title{Rigid objects in higher cluster categories}
\author{Anette Wrålsen}

\begin{abstract} 

We study maximal $m$-rigid objects in the $m$-cluster category $\mathcal C_H^m$ associated with a finite dimensional hereditary algebra $H$ with $n$ nonisomorphic simple modules. We show that all maximal $m$-rigid objects in these categories have exactly $n$ nonisomorphic indecomposable summands, and that any almost complete $m$-rigid object in $\mathcal C_H^m$ has exactly $m+1$ nonisomorphic complements. We also show that the maximal $m$-rigid objects and the $m$-cluster tilting objects in these categories coincide, and that the class of finite dimensional algebras associated with maximal $m$-rigid objects is closed under certain factor algebras.

\end{abstract}

\maketitle

\sloppy
%%%%%%%%%%%%%%%%%%%%%%%%%%%%%%%%%%%%%%%%%%%%%%%%%%%%%%%%%%%%%%%%%%%%%%%%%%%%%%%
\section*{Introduction}
%%%%%%%%%%%%%%%%%%%%%%%%%%%%%%%%%%%%%%%%%%%%%%%%%%%%%%%%%%%%%%%%%%%%%%%%%%%%%%%
Let $H$ be a finite dimensional hereditary algebra over an algebraically closed field $k$, such that $H$ has $n$ isoclasses of simple modules. Consider the bounded derived category $D^b(H)=\mathcal D$. Let $G=\tau^{-1}_{\mathcal D}[m]$ for some $m \geq 1$, where $\tau_{\mathcal D}$ is the AR-translation in $\mathcal D$ induced by the AR-translation of $\modu H$, and $[m]$ is the shift functor performed $m$ times. We can then form the orbit category $D^b(H)/G$ which we refer to as \emph{the $m$-cluster category} $\mathcal C^m_H$. All such categories $\mathcal C_H^m$ were shown to be triangulated in \cite{Keller}. The \emph{cluster category} $\mathcal C_H = D^b(H)/ \tau^{-1}[1]$ introduced and investigated in \cite{clustercat} is such a category, and the aim of this paper is to show that several important properties of the cluster category $\mathcal C_H$ generalize to $\mathcal C^m_H$.

We will study the \emph{maximal $m$-rigid objects} and \emph{$m$-cluster tilting objects} in these categories. An $m$-rigid object is an object $\hat T$ which is the direct sum of non-isomorphic indecomposable objects $\hat{T_1}, \hat{T_2},\ldots, \hat{T_r}$ such that $\Ext^i_{\mathcal C^m_H}(\hat{T_k}, \hat{T_l})=0$ for all $1 \leq i \leq m$, $1 \leq k,l \leq r$. For $1$-rigid objects we just write rigid. If $\hat T$ is maximal with this property, we call $\hat T$ a maximal $m$-rigid object (see \cite{Thomas}). An $m$-cluster tilting object is an $m$-rigid object which has the property that if $\Ext^i_{\mathcal C^m_H}(T,X)=0=\Ext^i_{\mathcal C^m_H}(X,T)$ for all $1\leq i\leq m$, then $X$ must be a summand of $T$ (see \cite{KR2}, \cite{I}). We will show that these classes of objects coincide, as was shown for $m=1$ in \cite{clustercat}.

In \cite{clustercat} it was also shown that in the cluster category all maximal rigid objects are induced by tilting modules over some hereditary algebra derived equivalent to $H$. In our setting we show that all maximal $m$-rigid objects in $\mathcal C^m_H$ are induced by exactly $n$ indecomposable nonisomorphic objects which constitute an $m$-rigid object in the subcategory $\modu H_0 \vee (\modu H_0)[1] \vee \ldots \vee (\modu H_0)[m-1]$ of $D^b(H_0)$, where $H_0$ is a hereditary algebra that is derived equivalent to $H$. This $m$-rigid object is maximal in the domain $\modu H_0 \vee (\modu H_0)[1] \vee \ldots \vee (\modu H_0)[m-1]\vee H_0[m]$. Tilting modules over $H$ will in particular induce maximal $m$-rigid objects. We also show that any $m$-rigid object $\hat X$ in $\mathcal C_H^m$ having $n-1$ nonisomorphic indecomposable summands has exactly $m+1$ \emph{complements}, i.e. nonisomorphic indecomposable objects $\hat Y$ such that $\hat X \coprod \hat Y$ is a maximal $m$-rigid object. This generalizes the property that any such object in $\mathcal C_H$ is known to have exactly $2$ complements (\cite{clustercat}, Theorem 5.1).

When $T$ is a maximal rigid object in $\mathcal C_H$, the associated algebra $\Gamma = \End_{\mathcal C_H}(T)^{op}$ is called a cluster-tilted algebra (see \cite{BMR1}). An important property of these algebras is that if $\Gamma$ is a cluster tilted algebra, so is $\Gamma/\Gamma e\Gamma$, where $e$ is an idempotent associated with a vertex for $\Gamma$. We also generalize this property to the endomorphism algebra of the maximal $m$-rigid objects in $\mathcal C_H^m$.

The main idea of this paper is to use the techniques of \cite{BMR} and generalized versions of these to be able to prove most of the main results by induction. If $\bar T = M\coprod T$ is a maximal $m$-rigid object and $M$ is indecomposable, we will look at the image of $T$ that results when we localise $\mathcal D$ with respect to the category $\mathcal M = \add \{M[i]|i\in \mathbb Z\}$ and see that the properties we are interested in transfer to the new category. Most of our results can be considered generalizations of analogous results for the cluster category, in particular found in \cite{clustercat} and \cite{BMR}.

In section 1 we look at how the maximal $m$-rigid objects in $\mathcal C_H^m$ are induced by corresponding $m$-rigid objects lying in $\mathcal D_G = \modu H \vee (\modu H)[1] \vee \ldots \vee (\modu H)[m-1] \vee H[m]$ in $\mathcal D$, a natural fundamental domain of the functor $G$. These objects are maximal $m$-rigid in this domain. We show that the maximal $m$-rigid objects in $\mathcal C_H^m$ can actually be assumed to be induced by $m$-rigid objects in $\mathcal D$ lying in the domain $\mathcal D_G^- = \modu H \vee (\modu H)[1] \vee \ldots \vee (\modu H)[m-1]$. This will be needed for the results in section 2. There we define and discuss the localisation of $\mathcal D$ with respect to an indecomposable object $M$, leading to the new category $\mathcal D_{\mathcal M}$ which is associated with a hereditary algebra $H'$. In particular we show that $m$-rigidity of an object in $\mathcal D_G$ transfers to the corresponding object in the localised category when we localise with respect to one of its indecomposable summands. Furthermore the image of an $m$-rigid object in $\mathcal D_G^-$ will actually be found in $\modu H' \vee (\modu H')[1] \vee \ldots \vee (\modu H')[m-1] \vee H'[m]$. We also show that the number of nonisomorphic indecomposable summands of the image of an $m$-rigid object in $\mathcal D_{\mathcal M}$ when we localise in one of the indecomposable summands is one less than the number of indecomposable summands in the original object. Section 3 contains the main results of the paper previously mentioned.

This work, which is part of the author's PhD thesis, started at the end of 2004 with first showing that a tilting $H$-module induces a maximal $m$-rigid object in $\mathcal C_H^m$. In the meantime many papers have appeared dealing with $m$-cluster categories and the more general class of $m$-Calabi Yau categories (see for example \cite{ABST}, \cite{BM1}, \cite{BM2}, \cite{HJ1}, \cite{HJ2}, \cite{IY}, \cite{KR1}, \cite{KR2}, \cite{Thomas}, \cite{Z}). In particular, the results on the number of indecomposable summands of maximal $m$-rigid objects and complements of almost complete such objects being respectively $n$ and $m+1$ when $H$ is of finite representation type was first shown in \cite{Thomas}. It has been proved in \cite{Z} for an arbitrary $H$ that the number of nonisomorphic indecomposable summands of $m$-cluster tilting objects is $n$. Since the maximal $m$-rigid objects and $m$-cluster tilting objects coincide, our result that the maximal $m$-rigid objects have $n$ summands gives a different approach to this result. The fact that maximal $m$-rigid objects and $m$-cluster tilting objects coincide and that almost complete $m$-cluster tilting objects have $m+1$ complements has been proved independently in \cite{ZZ} (with one inequality in \cite{IY}\cite{Z}).

I wish to thank my advisor Idun Reiten and Aslak Bakke Buan for many fruitful discussions and helpful comments on this paper. I would also like to thank the referee for many helpful suggestions.

%%%%%%%%%%%%%%%%%%%%%%%%%%%%%%%%%%%%%%%%%%%%%%%%%%%%%%%%%%%%%%%%%%%%%%%%%%%%%%%%%%%%%%%%%%%%%%%
\section{$m$-rigid objects in $\mathcal D$ inducing maximal $m$-rigid objects in $\mathcal C_H^m$}
%%%%%%%%%%%%%%%%%%%%%%%%%%%%%%%%%%%%%%%%%%%%%%%%%%%%%%%%%%%%%%%%%%%%%%%%%%%%%%%%%%%%%%%%%%%%%%%
In this section we will give some results linking maximal $m$-rigid objects in $\mathcal C_H^m$ with maximal $m$-rigid objects in $\mathcal D$ lying in $\mathcal D_G$ and $\mathcal D_G^-$. By maximal we here mean that they are maximal within the domain $\mathcal D_G$. We will also make an observation concerning a specific class of maximal $m$-rigid objects in $\mathcal C_H^m$, namely those induced by tilting modules over $H$.

The main result of this section is that any maximal $m$-rigid object in $\mathcal C_H^m$ is actually induced by a maximal $m$-rigid object in the smaller domain $\mathcal D_G^-$. To achieve this we may have to replace $H$ with a derived equivalent algebra $H_0$. This corresponds to the cluster category property that every cluster tilting object is induced by a tilting module over $H_0$ (and therefore has no summands in $H_0[1]$), and our proof is inspired by the proof of Theorem 3.3 in \cite{clustercat}.

Before we continue we recall the definition of a \emph{tilting module over a finite dimensional hereditary algebra $H$}. A module $T$ over $H$ is a tilting module if it is rigid and there is an exact sequence $0 \rightarrow H \rightarrow T^0 \rightarrow T^1 \rightarrow 0$ where $T^0$ and $T^1$ are in $\add T$. A module that is a direct summand in a tilting module is referred to as a \emph{partial tilting module}. There is also the dual notion of a \emph{cotilting module}, and it is shown in \cite{HR} that tilting modules and cotilting modules coincide for a hereditary algebra $H$.

There are several equivalent characterizations of tilting modules. One of them is that an $H$-module $T$ is a tilting module if it is maximal rigid in $\modu H$, another is that $T$ is a tilting module if and only if it is rigid and has $n$ nonisomorphic indecomposable direct summands, where $n$ is the number of non-isomorphic simple modules over $H$ (see \cite{HR}, Thm. 4.5). The general theory of tilting modules can be found in \cite{HR} and \cite{B}.

Given a maximal rigid object in the chosen fundamental domain of $\mathcal C_H^m$, any summand contained in a given shift will be induced by a partial tilting module. Therefore the object can have at most $n$ summands from each shift of the fundamental domain of $\mathcal C_H^m$, and since the fundamental domain of $\mathcal C_H^m$ is part of $m+1$ shifts in $\mathcal D$, the object can have at most $(m+1)n$ indecomposable non-isomorphic summands. Since $m$-rigid objects in particular are rigid, they also can have at most $(m+1)n$ indecomposable non-isomorphic summands.

\begin{definition}
We define the \emph{degree} of an object in $\mathcal D$ as follows: If $X \simeq Z[t]$ where $Z\in \modu H$ and $t \in \mathbb Z$, then the degree of $X$ is said to be $t$, written $\deg(X)=t$.
\end{definition}

The following lemma states that $m$-rigid objects in $\mathcal C^m_H$ are all induced by $m$-rigid objects in the chosen fundamental domain $\mathcal D_G$ in $\mathcal D$, and that the converse is also true.

\begin{lemma}\label{nullherogder}
Let $H$ be a hereditary algebra, and $X \in \mathcal D_G$ in $D^b(H)= \mathcal D$. Then $Hom_{\mathcal D}(X,X[i])=0$ for $i = 1, 2, \ldots, m$ if and only if $\Ext^i_{\mathcal C^m_H}(\widehat X,\widehat X)=0$ for $i=1, 2, \ldots, m$, where $\widehat X$ is the image of $X$ in $\mathcal C_H^m$.
\end{lemma}

\begin{proof}
\begin{align*}
\Ext^k_{\mathcal C^m_H}(\widehat{X},\widehat{X})	& \simeq \coprod_{t \in \mathbb{Z}} \Hom_{\mathcal{D}}(X,G^tX[k]) \\
            															& \simeq \Hom_{\mathcal{D}}(X,\tau_{\mathcal D}X[k-m])\coprod \Hom_{\mathcal D}(X,X[k]) \\
            															& \simeq D\Ext^1_{\mathcal{D}}(X[k-m],X)\coprod \Hom_{\mathcal D}(X,X[k]) \\
            															& \simeq D\Hom_{\mathcal{D}}(X[k-m],X[1])\coprod \Hom_{\mathcal D}(X,X[k]) \\
            															& \simeq D\Hom_{\mathcal{D}}(X,X[1+m-k])\coprod \Hom_{\mathcal D}(X,X[k])
\end{align*}

If we let $k$ run from $1$ to $m$ we see that the claim follows.
\end{proof}

Now we are ready for the main result of this section.

\begin{theorem}\label{erimskift}
Let $\widehat T$ be a (maximal) $m$-rigid object in $\mathcal C^m_H$. Then $\widehat T$ is induced by a (maximal) $m$-rigid object $T$ in $\modu H_0 \vee (\modu H_0)[1] \vee \ldots \vee (\modu H_0)[m-1]$ where $H_0$ is a hereditary algebra that is derived equivalent to $H$. Furthermore $T$ is also (maximal) $m$-rigid in $\modu H_0 \vee (\modu H_0)[1] \vee \ldots \vee (\modu H_0)[m-1]\vee H_0[m]$.
\end{theorem}

\begin{proof}
Let $\widehat T$ be a basic maximal $m$-rigid object in $\mathcal C^m_H$, induced by the indecomposable objects $T_1, T_2, \ldots, T_r$ in the chosen fundamental domain $\mathcal D_G$ in $\mathcal D$ of $\mathcal C_H^m$. First assume that $H$ is of infinite representation type. If no $T_l$ is a summand of $H[m]$, we are done. So assume that some $T_l$ is a summand of $H[m]$, and consider $\tau^{-1}_{\mathcal D} T$. If it has no summands in $H[m]$, we are also done since we can choose $H_0 \simeq H$, just with a different embedding into $\mathcal D$. Since we have already established that $T$ has at most $(m+1)n$ nonisomorphic indecomposable summands, $\tau^{-k}T$ will have no summands in $\modu H[m]$ for some $k$, $1 \leq k \leq (m+1)n$, and we can again choose $H_0$ isomorphic to $H$ but with a different embedding into $\mathcal D$. 

Next assume that $H$ is of finite representation type. Again, if no $T_k$ is a summand of $H[m]$ we are done. So assume that some $T_k$ is a summand of $H[m]$. If no $T_l$ is a summand of $H$, we can, as above, get the result by choosing $H_0 \simeq H$ but with a different embedding into $\mathcal D$.

Therefore assume that some $T_l$ is a summand of $H$. We first note that if all the objects induced by simple projective modules over $H$ are summands of $T$, no $T_j$ can be a summand of $\tau_{\mathcal D}^{-1}H$. This follows from the fact that if $P$ is a projective $H$-module, then $\Ext^1_{\mathcal D}(\tau_{\mathcal D}^{-1}P,S)\simeq \Hom_H(S,P)\neq 0$ for at least one simple projective $H$-module $S$. We want to show that if not all the projective simples of $\modu H$ are already summands of $T$, we can replace $H$ by a derived equivalent algebra $H_0$ such that this will be the case. 

Assume that $S$ is a simple projective object in $\modu H$ such that $S\notin \add T$. Then we claim that there is some path from $S$ to a summand $T_j$ of $T$, and hence a path in the AR-quiver of $\mathcal D$.

Since $S$ is not in $\add T$ and $S$ is $m$-rigid, there is some $i$, $1 \leq i \leq m$, such that $\Ext^i_{\mathcal C^m_H}(\widehat T, \widehat S)\not= 0$. We have
    
\begin{align*}
\Ext^i_{\mathcal C^m_H}(\widehat{T},\widehat{S}) 	& \simeq \Ext^1_{\mathcal C^m_H}(\widehat{T},\widehat{S[i-1]}) \\
                    											& \simeq D\Hom_{\mathcal C^m_H}(\widehat{S[i-1]},\tau_{\mathcal C^i_H}\widehat{T}) \\
                    											& \simeq \coprod_{t \in \mathbb{Z}} D\Hom_{\mathcal{D}}(S[i-1],G^t\tau_{\mathcal D}T) \\
                    											& \simeq \coprod_{t \in \mathbb{Z}} D\Hom_{\mathcal{D}}(S,G^t\tau_{\mathcal D}T[-i+1]).
\end{align*}

The degree of the indecomposable summands of $T$ is at most $m$, and so the degree of the indecomposable summands of $\tau_{\mathcal D}T[-i+1]$ will be at most $m-1$ since $i\geq 1$ and any summand of $T$ of degree $m$ will be projective. Therefore $\Hom_{\mathcal{D}}(S,G^t\tau_{\mathcal D}T[-i+1])$ must be zero when $t<0$. Similarly the degree of the indecomposable summands of $\tau_{\mathcal D}T[-i+1]$ is greater than or equal to $-m$ since the indecomposable summands of $T$ are of degree at least $0$, so $\Hom_{\mathcal{D}}(S,G^t\tau_{\mathcal D}T[-i+1])$ must be zero when $t\geq 2$. This means that at least one of $\Hom_{\mathcal{D}}(S,\tau_{\mathcal D}T[-i+1])$ and $\Hom_{\mathcal{D}}(S,T[m-i+1])$ is nonzero. But it cannot be $\Hom_{\mathcal{D}}(S,T[m-i+1])$ since $m-i+1 \geq 1$, and so no summand of $T[m-i+1]$ can be in the same shift as the simple projective $S$. Therefore we know that we have a path from $S$ to $\tau_{\mathcal D}T_j[-i+1]$ for some indecomposable summand $T_j$ of $T$, and consequently a path from $S$ to $T_j[-i+1]$. But we need a path from $S$ to $T_j$. If $X$ is an indecomposable object in $\modu H$, we know that there is a path from $X$ to $X[1]$ in $\mathcal D$, due to the following argument:

Let $S'$ be a simple module in $\modu H$ such that $S'$ is a composition factor of $X$. Then we know that $\Hom_H(X,I(S'))\not=0$ and $\Hom_H(P(S'),X)\not= 0$, where $I(S')$ is the injective envelope of $S'$ and $P(S')$ is its projective cover. Furthermore we know that $\tau_{\mathcal D} P(S')=I(S')[-1]$ in $\mathcal D$, and so we have a path $X \rightarrow I(S') \rightarrow P(S')[1] \rightarrow X[1]$. 

This means that there is a path from $T_j[t]$ to $T_j[t+1]$ for all $t\in \mathbb Z$, and so there must be a path from $T_j[-i+1]$ to $T_j$ giving us a path from $S$ to $T_j$.

Next let $\alpha (H)$ be the sum of the length of all paths in the AR-quiver of $\mathcal D$ from the simple projectives of $\modu H$ to some indecomposable summand of $T$. Since we have assumed that there is at least one simple projective not in $\add T$, we know that $\alpha (H)> 0$. Now, by if necessary replacing $H$ with a derived equivalent hereditary algebra, we can assume that $\alpha (H)$ is smallest possible. If there still is a simple projective object $S$ not in $\add T$, we can perform an APR-tilt (see \cite{APR} and \cite{BMR1}) using the basic tilting module $\tau^{-1}S\coprod P$ (where $H = S \coprod P$). This gives us the new algebra $\End_H(\tau^{-1}S\coprod P)^{\op}$ which is derived equivalent to $H$. Then $\alpha (\End_H(\tau^{-1}S\coprod P)^{\op})<\alpha (H)$, contradicting the minimality of $\alpha (H)$, and so all simple projectives must be in $\add T$. Now we choose $H_0$ derived equivalent to $H$ such that 
\begin{align*}
\tau_{\mathcal D}^{-2}(\modu H \vee (\modu H)[1] \vee \ldots \vee (\modu H)[m-1] \vee H[m])\\
 = \modu H_0 \vee (\modu H_0)[1] \vee \ldots \vee (\modu H_0)[m-1] \vee H_0[m].
\end{align*} 
Regarding $T$ as an object in $\modu H_0 \vee (\modu H_0)[1] \vee \ldots \vee (\modu H_0)[m-1] \vee H_0[m]$ we get the desired result.

We see that $\widehat T$ now will correspond to a (maximal) $m$-rigid object in $\mathcal C_{H_0}^m$, and so by Lemma \ref{nullherogder} $T$ is (maximal) $m$-rigid also in $\modu H_0 \vee (\modu H_0)[1] \vee \ldots \vee (\modu H_0)[m-1] \vee H_0[m]\vee H_0[m]$.

\end{proof}

We end this section by showing that tilting modules over $H$ induce maximal $m$-rigid objects in $\mathcal C_H^m$. Note that this was generalized in section 5.6 of \cite{KR1}. First we give an observation that will simplify some of the calculations.

\begin{lemma}\label{maksen}
Let $X$ and $Y$ be objects in $\mathcal D_G$ in $\mathcal D$. Then	$\Hom_{\mathcal D}(X,G^iY) = 0$ for all $i \neq 0,1$, and $\Hom_{\mathcal D}(X,G^jY)\neq 0$ for at most one of these values of $j$ if $m\geq 2$.
\end{lemma}

\begin{proof}
$\Hom_{\mathcal D}(X,G^iY) = \Hom_{\mathcal D}(X,\tau^{-i}Y[mi])$. When $i \leq -2$ this is obviously $0$ since then $\deg(G^{i}Y)<0$. When $i = -1$ we get the same unless $Y \in H[m]$. Then $\deg(Y[-m])=0$, but since $Y[-m]$ corresponds to a projective $H$-module, $\deg(G^{-1}Y)=\deg(\tau Y[-m])=-1$.

When $i \geq 2$ we see that even if $\deg(X)=m$ and $\deg(Y)=0$, we have $\deg(G^iY)>m+1$ and so $\Hom_{\mathcal D}(X,G^iY)$ will be zero.

Next, assume that both $\Hom_{\mathcal D}(X,Y)$ and $\Hom_{\mathcal D}(X,GY)$ are nonzero. We have 
\begin{align*}
\Hom_{\mathcal D}(X,GY)	& \simeq \Hom_{\mathcal D}(X,\tau^{-1}Y[m])	\\
												& \simeq \Hom_{\mathcal D}(\tau^2 X[-m],\tau Y)\\
										 		& \simeq D \Ext^1_{\mathcal D}(Y, \tau^2 X[-m])\\
										 		& \simeq D\Hom_{\mathcal D}(Y,\tau^2 X[-m+1]).
\end{align*} 

$\Hom_{\mathcal D}(X,Y) \neq 0$ implies that $\deg(X)\leq \deg(Y)$, and $\deg(\tau^2 X[-m+1])<\deg(X)$ when $m \geq 2$, so there can be no nonzero morphism from $Y$ to $\tau^2 X[-m+1]$, which gives us a contradiction.
\end{proof}

\begin{proposition}
Let $\widehat T$ be induced by a tilting module $T$ over $H$. Then $\widehat T$ is a maximal $m$-rigid object in $\mathcal C^m_H$.
\end{proposition}

\begin{proof}
Let $T \simeq T_1 \coprod T_2 \coprod \ldots \coprod T_n$ be a tilting module, where all the $T_i$ are indecomposable and	nonisomorphic. We first check that $\Ext^i_{\mathcal{D}/G}(\widehat{T},\widehat{T})=0$ for $1\leq i \leq m$, which is the same as checking that $\Ext^i_{\mathcal{D}/G}(\widehat{T_k},\widehat{T_l})=0$ for any $1\leq k,l \leq n$:

\begin{align*}
\Ext^i_{\mathcal{D}/G}(\widehat{T_k},\widehat{T_l})	& \simeq \coprod_{t \in \mathbb{Z}} \Ext^i_{\mathcal{D}}(T_k,G^tT_l) \\
            																& \simeq \coprod_{t \in \mathbb{Z}} \Hom_{\mathcal{D}}(T_k,G^tT_l[i]) \\
            																& \simeq \Hom_{\mathcal{D}}(T_k,T_l[i]) \coprod \Hom_{\mathcal{D}}(T_k,\tau^{-1}T_l[i+m])
\end{align*}

The first summand corresponds to $\Ext^i_H(T_k,T_l)$ for $1 \leq i \leq m$ which is zero by assumption, since $T$ is a tilting module and $H$ is hereditary. Since $i+1 \geq 2$, we get that $\deg(\tau^{-1}T_l[i+m])\geq 2$ and so the second summand is zero.

To show that $\widehat T$ is maximal, we consider all indecomposable rigid objects $X$ of $\mathcal D_G$ not in $\{T_1, T_2, \dots T_n\}$ and show that none of them can be added to $\widehat T$ whilst keeping it an 	$m$-rigid object.

First assume that $\widehat X \simeq \widehat{T_k[i]}$ for some $k$, and $1\leq i \leq m$. We then see that    
\begin{align*}
\Ext^i_{\mathcal D/G}(\widehat{X}, \widehat{T_k}) 	& \simeq \coprod_{t \in \mathbb Z} \Ext^i_{\mathcal D}(T_k[i], G^tT_k) \\
            																				& \simeq \coprod_{t \in \mathbb Z} \Hom_{\mathcal D}(T_k[i], G^tT_k[i]) \\
            																				& \simeq \Hom_{\mathcal D}(T_k[i], T_k[i])\\
           																					& \simeq \Hom_H(T_k, T_k)\neq 0.
\end{align*}

Next assume that $\widehat X \simeq \widehat{Z[i]}$, $Z \in \modu H$, $0 \leq i \leq m$, where $Z$ is \emph{not} isomorphic to any summand of $T$. If $\deg \widehat X = 0$ obviously $\widehat X$ cannot extend $\widehat T$ (otherwise $T$ would not be a tilting module over $H$), so we can assume that $i \geq 1$. We now know that there exists some $l$, $1 \leq l \leq n$, such that $\Ext^1_H(Z,T_l)$ or $\Hom_H(Z,T_l)$ is nonzero since $T$ is also a cotilting module (this follows from for instance \cite{HR}). If $\Ext^1_H(Z,T_l)$ is nonzero, we see that 
\begin{align*}
\Ext^{i+1}_{\mathcal D/G}(\widehat{X}, \widehat{T_l}) 	& \simeq \coprod_{t \in \mathbb Z} \Ext^{i+1}_{\mathcal D}(Z[i], G^tT_l) \\
            																						& \simeq \Ext^1_{\mathcal D}(Z[i], T_l[i]) \\
            																						& \simeq \Ext^1_H(Z, T_l).
\end{align*}

If $i=m$, then $\widehat X \simeq \widehat{Z[m]}$ such that $Z$ is a projective $H$-module, and so $\Ext^1_H(Z,Y)=0$ for any $H$-module $Y$.

We also know by Lemma \ref{maksen} that $\Ext^i_{\mathcal D/G}(\widehat{X},\widehat{T_k}) \simeq \coprod_{t \in \mathbb Z} \Ext^i_{\mathcal D}(G^tZ[i],T_l) \simeq \Hom_H(X, T_l)$. This concludes the proof.

\end{proof}

%%%%%%%%%%%%%%%%%%%%%%%%%%%%%%%%%%%%%%%%%%%%%%%%%%%%%%%%%%%%%%%%%%%%%%%%%%%%
\section{Localising with respect to a rigid object}
%%%%%%%%%%%%%%%%%%%%%%%%%%%%%%%%%%%%%%%%%%%%%%%%%%%%%%%%%%%%%%%%%%%%%%%%%%%%

Let $M$ be an indecomposable rigid $H$-module for some finite dimensional hereditary algebra $H$. Then it is known that the full subcategory of $\modu H$	consisting of the $H$-modules $X$ such that $\Hom_H(M,X)=0=\Ext^1_H(M,X)$, which we will refer to as $\mathcal U_M$, is an exact subcategory of $\modu H$. This subcategory is equivalent to $\modu H'$, where $H'$ is the endomorphism ring of the direct sum of the nonisomorphic indecomposable projective objects of $\mathcal U_M$. $H'$ is a hereditary algebra with one fewer isoclass of simple modules than $H$. Details on this can be found for instance in \cite{Happel}.

We will now describe how \emph{localisation of $\mathcal D$ with respect to an indecomposable rigid object $M$} gives rise to the derived category of a hereditary algebra $H'$ with one fewer isoclass of simple modules. We will only give the definitions needed for our purposes; more general notions and details on this construction can be found in \cite{BMR}. Combining this with Theorem \ref{erimskift}, we get the mathematical machinery needed to make our inductive construction. 

So let $M$ be as above, and let $\mathcal M=\add\{M[i]|i\in\mathbb Z\}$, a thick triangulated subcategory of $\mathcal D$. Then we can  localise $\mathcal D$ at $\mathcal M$ to form a new category $\mathcal D_{\mathcal M}$. There is a localisation functor $L_{\mathcal M}:\mathcal D \rightarrow \mathcal D_{\mathcal M}$, which is a canonical exact triangle functor with the property that $L_{\mathcal M}(M') = 0$ for any $M' \in \mathcal M$, and it is universal with respect to this property. We will use the notation $L_{\mathcal M}(X)=\widetilde X$ for $X \in \mathcal D$.

It is shown in \cite{BMR} that there is an exact equivalence between $\mathcal D_{\mathcal M}$ and $\mathcal D_0$, where $\mathcal D_0=\add\{U|U\in\mathcal D, \Hom_{\mathcal D}(M,U[i])=0 \mbox{ for all } i\in \mathbb Z\}$. Note that if $M$ is an $H$-module, then $\mathcal D_0=\add\{U[i]| U\in \mathcal U_M, i\in \mathbb Z\}$. The following proposition gives an important consequence of this equivalence (this is an adjustment of Proposition 2.2 in \cite{BMR} for our purposes, see also \cite{V1}, Ch. 2, 5-3, \cite{V2}):

\begin{proposition}\label{dnulliso}
Let $\mathcal D$, $\mathcal M$, $\mathcal D_0$ and $\mathcal D_{\mathcal M}$ be as above. Given an object $Y$ in $\mathcal D$, then $Y\in\mathcal D_0$ if and only if for every object $X$ of $\mathcal D$ the canonical map $\Hom_{\mathcal D}(X,Y) \rightarrow \Hom_{\mathcal D_{\mathcal M}}(L_{\mathcal M}(X),L_{\mathcal M}(Y))$ is an isomorphism.
\end{proposition}

Since $L_{\mathcal M}$ is a triangle functor, we have $L_{\mathcal M}(X[i])\simeq L_{\mathcal M}(X)[i]$ in $\mathcal D_{\mathcal M}$ for any $X \in \mathcal D$. We will also need Lemma 2.14 of \cite{BMR}:

\begin{lemma}\label{dnullogtau}
Let $X$ be an indecomposable object in $\mathcal D_0 \subset \mathcal D$. Then $\widetilde X$ is indecomposable and $\widetilde{\tau^{-1}_{\mathcal D}X} \simeq \tau^{-1}_{\mathcal D'}\widetilde X$.
\end{lemma}

The following lemma gives a very fundamental connection between maps in $D_{\mathcal M}$ and $\mathcal D$ (see \cite{V1}, \cite{V2}). 

\begin{lemma}\label{nullilokalisert}
For any map $f$ in $\mathcal D$, $L_{\mathcal M}(f) = 0$ if and only if $f$ factors through an object in $\mathcal M$.
\end{lemma}

In \cite{BMR} we find the following result about $\mathcal D_{\mathcal M}$.

\begin{theorem}
Let $H$ be a hereditary algebra with $n$ simple modules up to isomorphism. Let $M$ be an indecomposable $H$-module with $\Ext^1_H(M,M)=0$, and let $\mathcal M$ be as above. Then $\mathcal D_{\mathcal M}$ is equivalent to the derived category of a hereditary algebra with $n-1$ simple modules (up to isomorphism).
\end{theorem}

We now give a technical lemma that will simplify some of the arguments in this section.

\begin{lemma}\label{blirnull}
Let $M \coprod X$ be an $m$-rigid object in $\mathcal D_G$ such that $M$ is indecomposable, and let $M_X \rightarrow X$ be the minimal right $\mathcal M$-approximation of $X$ in $\mathcal D$ inducing the triangle 
\[M_X \xrightarrow{f} X \xrightarrow{g} Y \rightarrow\]
Then $Y \in \mathcal D_0$, and $\Hom_{\mathcal D}(X,M_X[t])=0$ for $t\geq 1$.
\end{lemma}

\begin{proof}
Since $f$ is right minimal, $X \in \mathcal D_G$ and $\Hom_{\mathcal D}(M,X[k])\simeq \Hom_{\mathcal D}(M[-k],X)=0$ for $1 \leq k \leq m$, then $M_X$ must be in $\add\{M[j]|j = 0, 1, \ldots, m\}$.

One way to show that $Y \in \mathcal D_0$ is to show that $\Hom_{\mathcal D}(M,Y[j]) = 0$ for all $j$. Applying the functor $\Hom_{\mathcal D}(M,-)$ to the above triangle, we get the long exact sequence
\begin{align*}
\Hom_{\mathcal D}(M,M_X[j]) \xrightarrow{f_*} \Hom_{\mathcal D}(M,X[j]) & \xrightarrow{g_*} \Hom_{\mathcal D}(M,Y[j]) \\
																																				&	\xrightarrow{h_*} \Hom_{\mathcal D}(M,M_X[j+1])            
\end{align*}

Since $f$ is a right $\mathcal M$-approximation, $f_*$ is an epimorphism. Furthermore, since $M$ is rigid and $f$ is right minimal, it must be a monomorphism (since any morphism from $M$ to itself is either zero or an isomorphism). This means that $g_*$ and $h_*$ always will be zero, forcing $\Hom_{\mathcal D}(M,Y[j])$ to be zero for all $j$. Therefore $Y \in \mathcal D_0$.

Next consider $\Hom_{\mathcal D}(X,M_X[t])=0$ when $t\geq 1$. Since $\Hom_{\mathcal D}(X,M[k])=0$ for $k=1, 2, \ldots, m$ and $M_X \in \add\{M[j]|j = 0, 1, \ldots, m\}$, the indecomposable summands of $M_X[t]$ will be of the form $M[j+t]$ where $j+t\geq 1$. By assumption $\Hom_{\mathcal D}(X,M[j+t])=0$ when $1 \leq j+1 \leq m$. If $j+t = m+1$, then $\Hom_{\mathcal D}(X,M[j+t])$ can only be nonzero if $X$ is of degree $m$. But if $X$ is of degree $m$ it is of the form $P[m]$ where $P$ is a projective $H$-module, and so $\Hom_{\mathcal D}(X,M[j+t])\simeq \Hom_{\mathcal D}(P[m],M[m+1]) \simeq \Ext^1_H(P,M) = 0$. Since $X \in \mathcal D_G$, we have $\Hom_{\mathcal D}(X,M[j+t])=0$ when $j+t>m+1$, and so we see that $\Hom_{\mathcal D}(X,M_X[t])=0$ for all $t\geq 1$.

\end{proof}

We will now investigate what happens to an $m$-rigid object in $\mathcal D$ when we localise with respect to one of its indecomposable summands. By Theorem $\ref{erimskift}$ it is enough to consider $m$-rigid objects in $\mathcal D_G^-$. The next result is an analogue of Lemma 2.10 in \cite{BMR}. First we define what we mean by a \emph{complement} of an $m$-rigid object in our setting.

\begin{definition}
Let $M$ be an $m$-rigid object in $\mathcal D_G$. We call $T$ \emph{a complement of $M$ in $\mathcal D_G$} if no summand of $T$ is in $\add M$, and $M\coprod T$ is a maximal $m$-rigid object in $\mathcal D_G$.
\end{definition}

\begin{lemma}\label{tmerketeri}
Let $M$ be a rigid indecomposable object in $\mathcal D_G^-$, and let $T$ be a complement of $M$ in $\mathcal D_G$ such that $T \in \mathcal D_G^-$. Let $D_{\mathcal M} = \mathcal L_{\mathcal M}(\mathcal D)$, and $H'$ be the hereditary algebra corresponding to $\mathcal U_X$ where $X$ is the $H$-module such that $M = X[k]$, $0 \leq k \leq m-1$, as previously described. Then we have the following:
\begin{itemize}
\item[(a)]   
$L_{\mathcal M}(T)=\widetilde T$ is in $\modu H' \vee (\modu H')[1] \vee \ldots \vee (\modu H')[m-1] \vee H'[m]$.
    
\item[(b)]   
$\Hom_{\mathcal D_{\mathcal M}}(\widetilde T, \widetilde T[k])= 0$ for $1 \leq k \leq m$. 
\end{itemize}

\end{lemma}

\begin{proof}

(a) Let $f: M' \rightarrow T$ be a minimal right $\mathcal M$-approximation of $T$ in $\mathcal D$ inducing the triangle 
\[M'\xrightarrow{f} T \xrightarrow{g} U_T\rightarrow.\] 
We know by Lemma \ref{blirnull} that $U_T\in \mathcal D_0$. We can also conclude that all indecomposable summands of $U_T$ must have a degree in $\{0,1,\ldots,m\}$ since $f$ is a minimal right $\mathcal M$-approximation.

Since $T \simeq U_T$ in $\mathcal D_{\mathcal M}$, we know that if $U_T$ has a nonzero summand $U_m$ in $(\modu H)[m]$, $T$ must have a nonzero summand $T_{m-1}$ in $(\modu H)[m-1]$ such that there exists a triangle $M''\rightarrow T_{m-1}\rightarrow U_m \rightarrow$ with $M'' \in \mathcal M$. Since $T_{m-1}\simeq T^*[m-1]$ for some $T^* \in \modu H$, we know that $M''\in \add\{M[m-2],M[m-1]\}$. But since $U_m\simeq U^*[m]$ for some $U^*\in \modu H$ and therefore $\Hom_{\mathcal D}(U_m,M[m-2][1])=0$, we see that $M''\simeq M^*[m-1]$ for some $M^*\in \add M$ since $M$ is not a summand of $T$.

This means that the triangle $M''\rightarrow T_{m-1}\rightarrow U_m \rightarrow$ is induced by a map $M^*\xrightarrow{\alpha} T^*$ in $\modu H$. This again means that $U^*\simeq \Coker \alpha [-1] \coprod \Ker \alpha$, and since $\deg(U_m)=m$ we can conclude that $U_m \simeq U^*[m] \simeq (\Ker \alpha)[m]$. 

Since $\Ker \alpha$ is a submodule of $M^*$, the inclusion $\Ker \alpha \hookrightarrow M$ induces the exact sequence $\Ext^1_H(M,X) \rightarrow \Ext^1_H(\Ker \alpha,X) \rightarrow 0$ for any $X \in \modu H$. In particular this means that $\Ext^1_H(\Ker \alpha,X) = 0$ for all $X \in \mathcal U_M$, and so $\Ker \alpha$ must be projective in $\mathcal U_M$. Therefore $U_m \simeq P[m]$ for some projective object $P$ in $\mathcal U_M$, and so we see that $U_T$ must be in $\modu H' \vee (\modu H')[1] \vee \ldots \vee (\modu H')[m-1] \vee H'[m]$.

\noindent (b) We can apply the functor $\Hom_{\mathcal D}(T,-)$ to the triangle from (a) to obtain the exact sequence
\[\Hom_{\mathcal D}(T,T[k]) \rightarrow \Hom_{\mathcal D}(T,U_T[k])\rightarrow\Hom_{\mathcal D}(T,M'[k+1]) \]
where $1\leq k \leq m$.

By Lemma \ref{blirnull}, we have that $\Hom_{\mathcal D}(T,M'[k+1])=0$ since $k\geq 1$. By assumption $\Hom_{\mathcal D}(T,T[k])=0$ for $k=1, 2, \ldots, m$, so $\Hom_{\mathcal D}(T,U_T[k])=0$ for all such $k$. Since $U_T[k]$ is in $\mathcal D_0$, Proposition \ref{dnulliso} assures us that $\Hom_{\mathcal D_{\mathcal M}}(\widetilde T,\widetilde U_T[k])=0$, and since $\widetilde T \simeq \widetilde{U_T}$ this means that $\Hom_{\mathcal D_{\mathcal M}}(\widetilde T, \widetilde T[k]) = 0$.

\end{proof}

Next we give a proposition that will prove to be very useful for demonstrating the maximality of the image of maximal $m$-rigid objects of $\mathcal D_G$ in $\mathcal C_{H}^m$.

\begin{proposition}\label{lagtriangel}

Let $M$ be an indecomposable rigid $H$-module, and let $Y\in \mathcal D_0$ be an indecomposable rigid object such that $\Hom_{\mathcal D}(Y,M[i])\neq 0$ for some $i\in \mathbb Z$. Then $\widetilde Y \simeq \widetilde X$ for some indecomposable rigid object $X \notin \mathcal D_0$ such that $\Hom_{\mathcal D}(M,X[t])=0$ for $t\neq 1-i$, and $\Hom_{\mathcal D}(X,M[t])=0$ for all $t$.

\end{proposition}

\begin{proof}
  
Let $Y \xrightarrow{f} M'$ be the minimal left $\mathcal M$-approximation of $Y$. Then $\deg(Y)$ is either $i$ or $i-1$, so $M'=M_1[i-1]\coprod M_2[i] \coprod M_3[i+1]$ for some $M_1, M_2, M_3 \in \add M$. Then $f$ can be completed to a triangle 
\begin{equation}\label{tri} X \rightarrow Y \xrightarrow{f} M'\rightarrow.\end{equation}

We have $\widetilde X \simeq \widetilde Y$, and $\widetilde Y$ will be rigid by Proposition \ref{dnulliso} since $Y\in \mathcal D_0$. Since $Y$ is indecomposable and $\widetilde X \simeq \widetilde Y$, we have that $X$ must be of the form $X^* \coprod M^*$ where $X^*$ is indecomposable and $M^* \in \mathcal M$. But $Y \in \mathcal D_0$, and so $\Hom_{\mathcal D}(M^*,Y)=0$. Since $f$ is a minimal left $\mathcal M$-approximation and the composition $Y \xrightarrow{f} M' \rightarrow (X^*\coprod M^*)[1]$ induced by (\ref{tri}) is zero, $M'\rightarrow M^*[1]$ must be zero (since $M$ is rigid and any nonzero map from $M$ to itself is an isomorphism). Therefore $M^*$ must be zero, and so $X$ is indecomposable.

Next we note that $\deg(X)$ must be either $i$ or $i-1$. This is because $f$ being a minimal left $\mathcal M$-approximation and $Y \ncong M[i]$ means that $\Hom_{\mathcal D}(M[i],X[1])$ is nonzero, so $X$ cannot be of degree less than $i-1$. It also follows from this that $X\notin \mathcal D_0$. Now only $\Hom_{\mathcal D}(X,M[i-1])$, $\Hom_{\mathcal D}(X,M[i])$ and $\Hom_{\mathcal D}(X,M[i+1])$ can possibly be nonzero and so we only have to check these. To do this we apply the functor $\Hom_{\mathcal D}(-,M)$ to (\ref{tri}). We then get the long exact sequence

\begin{align*}
	0 & \rightarrow \Hom_{\mathcal D}(M'[-i+1],M) \rightarrow \Hom_{\mathcal D}(Y[-i+1],M) \rightarrow \Hom_{\mathcal D}(X[-i+1],M)\\
		& \rightarrow \Hom_{\mathcal D}(M'[-i],M)   \rightarrow \Hom_{\mathcal D}(Y[-i],M)   \rightarrow \Hom_{\mathcal D}(X[-i],M)\\
    & \rightarrow \Hom_{\mathcal D}(M'[-i-1],M) \rightarrow \Hom_{\mathcal D}(Y[-i-1],M) \rightarrow \Hom_{\mathcal D}(X[-i-1],M)\\	
    & \rightarrow 0.     
\end{align*}

The map $\Hom_{\mathcal D}(M'[-i-1],M) \rightarrow \Hom_{\mathcal D}(Y[-i-1],M)$ is induced by $f$. Since $f$ is a left $\mathcal M$-approximation, it must be an epimorphism. But since $M$ is rigid and $f$ is left minimal, the map is a monomorphism too. The same reasoning goes for the maps $\Hom_{\mathcal D}(M'[-i],M) \rightarrow \Hom_{\mathcal D}(Y[-i],M)$ and $\Hom_{\mathcal D}(M'[-i+1],M) \rightarrow \Hom_{\mathcal D}(Y[-i+1],M)$, and so we can conclude that $\Hom_{\mathcal D}(X,M[i-1])$, $\Hom_{\mathcal D}(X,M[i])$ and $\Hom_{\mathcal D}(X,M[i+1])$ all must be zero since they are caught between zero maps in the long exact sequence. 

Assume that $X$ is not rigid. Then there is some nonzero map $g:X\rightarrow X[1]$. Since $\widetilde X$ is rigid in $\mathcal D_{\mathcal M}$, this map must factor through some object in $\mathcal M$, but we have just shown that $\Hom_{\mathcal D}(X,M'')=0$ for all $M''\in \mathcal M$, so this is impossible. 

Finally we consider $\Hom_{\mathcal D}(M,X[t])$ for $t\neq i-1$. If we apply the functor $\Hom_{\mathcal D}(M,-)$ to (\ref{tri}), we get the long exact sequence
\[ \Hom_{\mathcal D}(M,Y[t-1]) \rightarrow \Hom_{\mathcal D}(M,M'[t-1]) \rightarrow \Hom_{\mathcal D}(M,X[t]) \rightarrow \Hom_{\mathcal D}(M,Y[t]).\]
Since $Y \in \mathcal D_0$, we see that $\Hom_{\mathcal D}(M,M'[t-1])\simeq \Hom_{\mathcal D}(M,X[t])$ for all $t$. This means that $\Hom_{\mathcal D}(M,X[t])$ will be nonzero if and only if $M[1-t]$ is a summand in $M'$.

Now, since $\Hom_{\mathcal D}(X,Y)$ is nonzero and $\deg(X)$ and $\deg(Y)$ are either $i$ or $i-1$, we will have one of the two following cases. Either $\deg(Y)>\deg(X)$ or $\deg(Y)=\deg(X)$. In the first case $\deg(Y)=i$ and $\deg(X)=i-1$. Since $f$ is a minimal left approximation, this means that $M_1$ must be zero since $Y$ cannot map to anything of degree $i-1$. Similarly $X[1]$ is of degree $i$, and cannot be mapped to by anything of degree $i+1$. This means that $M_3$ also must be zero. Then $M'$ must be equal to $M_2[i]$, and only $\Hom_{\mathcal D}(M,X[1-i])$ can be nonzero.

If $\deg(Y)=\deg(X)$, we can apply the functor $\Hom_{\mathcal D}(-,X)$ to the triangle (\ref{tri}) to get the exact sequence
\[ \Hom_{\mathcal D}(X[-1],X) \rightarrow \Hom_{\mathcal D}(M'[-2],X) \rightarrow \Hom_{\mathcal D}(Y[-2],X).\]
Since $X$ is rigid and $\deg(Y)=\deg(X)$, we see that $\Hom_{\mathcal D}(M'[-2],X)=0$. We have that $M'[-2]=M_1[i-3]\coprod M_2[i-2] \coprod M_3[i-1]$, which means that in particular $\Hom_{\mathcal D}(M_2[i-2],X)=0$ and $\Hom_{\mathcal D}(M_3[i-1],X)=0$. Since $M_2$ is nonzero by assumption, the first identity means that $\Hom_{\mathcal D}(M[i-2],X)=0$, and so $M_1=0$. We have also already shown that $\Hom_{\mathcal D}(M[i-1],X)\neq 0$, so the second identity can only hold if $M_3=0$. Again we get that $M'=M_2[i]$, and $\Hom_{\mathcal D}(M,X[t])$ can only be nonzero for $t=1-i$. 
  
\end{proof}

Actually we see that this implies that $\Hom_{\mathcal D}(Y,M[t])$ can only be nonzero for $t=i$, and so $\Hom_{\mathcal D}(Y,M[t])$ can be nonzero for at most one $t\in \mathbb Z$ when $Y\in \mathcal D_0$.

The following result is a generalization of Proposition 2.12 of \cite{BMR}.

\begin{proposition}\label{liktantsummander}
Let $\bar T = M \coprod T$ be an $m$-rigid object in $\mathcal D_G$ in $\mathcal D$ such that $M$ is indecomposable and $M \notin \add T$. Then the image $\widehat T$ of $\widetilde T$ in $\mathcal C^m_{H'}$ is an $m$-rigid object with the same number of nonisomorphic indecomposable summands as $T$.
\end{proposition}

\begin{proof}

By Theorem \ref{erimskift} we can assume that $\bar T$ is in $\mathcal D_G^-$. By Lemma \ref{tmerketeri} we then know that $\widetilde T$ is in $\modu H' \vee (\modu H')[1] \ldots \vee (\modu H')[m-1] \vee H'[m]$, and Lemma \ref{nullherogder} assures us that $\Ext^k_{\mathcal C^m_{H'}}(\widehat T, \widehat T) = 0$ for $k = 1, 2, \ldots, m$. It remains to show that nonisomorphic indecomposable summands of $T$ are sent to nonisomorphic indecomposable objects in $\mathcal D_{\mathcal M}$.

Let $T_a$ be an indecomposable summand of $T$. If $T_a \in \mathcal D_0$, we know by Proposition \ref{dnulliso} that $\Hom_{\mathcal D_{\mathcal M}}(\widetilde T_a, \widetilde T_a) \simeq \Hom_{\mathcal D}(T_a,T_a)$, and so $\widetilde T_a$ must be indecomposable. If $T_a$ is not in $\mathcal D_0$, we can take the minimal right $\mathcal M$-approximation of $T_a$ and extend it to a triangle 
\begin{equation}\label{a} M_a\rightarrow T_a \rightarrow Y_a \rightarrow \end{equation} 
in $\mathcal D$. Then $Y_a$ will be in $\mathcal D_0$ by Lemma \ref{blirnull}. Furthermore $Y_a$ will be nonzero, and so $\Hom_{\mathcal D_{\mathcal M}}(\widetilde T_a, \widetilde T_a)\simeq \Hom_{\mathcal D_{\mathcal M}}(\widetilde Y_a, \widetilde Y_a)\simeq \Hom_{\mathcal D}(Y_a, Y_a)$ will be nonzero. We also have that $Y_a$ will be indecomposable. This is because $\Hom_{\mathcal D}(T_a,Y_a)\simeq\Hom_{\mathcal D}(Y_a,Y_a)$ by Proposition \ref{dnulliso}, and $\Hom_{\mathcal D}(T_a,T_a)\simeq\Hom_{\mathcal D}(T_a,Y_a)$. To see the last isomorphism, apply $\Hom_{\mathcal D}(T_a,-)$ to (\ref{a}) to get the exact sequence
\[\Hom_{\mathcal D}(T_a,T_a)\rightarrow\Hom_{\mathcal D}(T_a,Y_a)\rightarrow\Hom_{\mathcal D}(T_a,M_a[1]).\]

Since the map $\Hom_{\mathcal D}(T_a,T_a)\rightarrow\Hom_{\mathcal D}(T_a,Y_a)$ is nonzero and $\dim \Hom_{\mathcal D}(T_a,T_a)=1$ (since $T_a$ is rigid), this must be a monomorphism. Since $\Hom_{\mathcal D}(T_a,M_a[1])=0$ by Lemma \ref{blirnull} we get that it is an isomorphism.

Now assume that $T_b$ is an indecomposable summand of $T$ which is not isomorphic to $T_a$. We need to show that $\widetilde T_a$ is not isomorphic to $\widetilde T_b$. First we demonstrate that every nonzero map in $\Hom_{\mathcal D_{\mathcal M}}(\widetilde T_a, \widetilde T_b)$ is induced by a nonzero map in $\Hom_{\mathcal D}(T_a,T_b)$. 

If $T_b \in \mathcal D_0$, this holds by Proposition \ref{dnulliso}. If not, we can use the minimal right $\mathcal M$-approximation of $T_b$ to get a triangle 
\begin{equation}\label{b}M_b\rightarrow T_b \rightarrow Y_b \rightarrow.\end{equation} 
By Lemma \ref{blirnull} we have $Y_b \in \mathcal D_0$. We have that $\Hom_{\mathcal D_{\mathcal M}}(\widetilde T_a, \widetilde T_b) \simeq \Hom_{\mathcal D_{\mathcal M}}(\widetilde T_a, \widetilde Y_b)\simeq \Hom_{\mathcal D}(T_a, Y_b)$, so there is a nonzero map from $T_a$ to $Y_b$ in $\mathcal D$.

If we now apply the functor $\Hom_{\mathcal D}(T_a,-)$ to (\ref{b}), we get the following exact sequence:
\[\Hom_{\mathcal D}(T_a,T_b) \xrightarrow{g^*} \Hom_{\mathcal D}(T_a,Y_b) \rightarrow \Hom_{\mathcal D}(T_a,M_b[1]) \]

By Lemma \ref{blirnull} we know that $\Hom_{\mathcal D}(T_a,M_b[1])=0$, and so $g^*$ will be an epimorphism. This means that we must have a nonzero map from $T_a$ to $T_b$ for every nonzero map from $\widetilde T_a$ to $\widetilde T_b$.

Finally we need to rule out the possibility that a map $T_a \xrightarrow{\alpha}T_b$ in $\mathcal D$ can correspond to an isomorphism in $\mathcal D_{\mathcal M}$. If it does, it induces a triangle $T_a \xrightarrow{\alpha} T_b \rightarrow M' \rightarrow T_a[1] \rightarrow$ in $\mathcal D$ with $M'$ in $\mathcal M$. But by Lemma \ref{blirnull} $\Hom_{\mathcal D}(T_b,M[k])=0$ for all $k\geq 1$, so $M'$ would have to be in $\add\{M[k]|k \leq 0\}$.

But the same lemma also tells us that $\Hom_{\mathcal D}(M,T_a[l])\simeq \Hom_{\mathcal D}(M[-l],T_a)=0$ for all $l\geq 1$ since $\deg(M)$ is assumed to be less than or equal to $m-1$. So $\Hom_{\mathcal D}(M[k],T_a[1])\simeq\Hom_{\mathcal D}(M[k-1],T_a)=0$ for all $k\leq 0$. In other words, for $\alpha$ to correspond to an isomorphism $M'$ must be zero, which would contradict $T_a$ and $T_b$ being nonisomorphic. By Lemma \ref{tmerketeri} we now get that $\widehat T$ is preserved as an $m$-rigid object in $\mathcal C^m_H$ with its original number of nonisomorphic indecomposable summands.
\end{proof}

%%%%%%%%%%%%%%%%%%%%%%%%%%%%%%%%%%%%%%%%%%%%%%%%%%%%%%%%%%%%%%%%%%%%%%%%%%%%%%%
\section{Main results}
%%%%%%%%%%%%%%%%%%%%%%%%%%%%%%%%%%%%%%%%%%%%%%%%%%%%%%%%%%%%%%%%%%%%%%%%%%%%%%%

In this section we will prove the main results of this paper. We start by recalling the setting. Let $T$ be an $m$-rigid object in $\mathcal C_H^m$. Then it corresponds to an object $\bar T \in \mathcal D_G$ in $\mathcal D$. We can localise $\mathcal D$ with respect to any indecomposable rigid summand $M$ of $T$, and the resulting category $\mathcal D_{\mathcal M}$ is the derived category of a hereditary algebra $H'$ that has one fewer isoclass of simple modules than $H$.

The first theorem states that maximal $m$-rigid objects in $\mathcal D_G$ induce maximal $m$-rigid objects in $\mathcal C_{H'}^m$.

\begin{theorem}\label{teorem1}
Let $\bar T = M \coprod T$ be a maximal $m$-rigid object in $\mathcal D_G$ in $\mathcal D$ such that $M$ is indecomposable and $M\notin \add T$. Then the image $\widehat T$ of $\widetilde T$ in $\mathcal C^m_{H'}$ is a maximal $m$-rigid object.
\end{theorem}

\begin{proof}
We know by Proposition \ref{liktantsummander} that $\widehat T$ is an $m$-rigid object in $\mathcal C_{H'}^m$, so we only need to show that it is maximal. Assume that it is not.  Then there is some indecomposable rigid object $\widehat C \in \mathcal C^m_{H'}$ such that $\widehat T \coprod \widehat C$ is $m$-rigid in $\mathcal C^m_{H'}$ and $\widehat C$ is not isomorphic to any summand of $\widehat T$. The object $\widehat C$ is induced by an object $\widetilde C$ in $\modu H' \vee (\modu H')[1] \ldots \vee (\modu H')[m-1] \vee H'[m]$ in $\mathcal D_{\mathcal M}$. This object can be lifted to an indecomposable object in $\mathcal D$, and in particular to some indecomposable object $C$ in $\mathcal D_0$. $C$ must be rigid by Proposition \ref{dnulliso}. The same proposition gives us that $\Hom_{\mathcal D}(\bar T,C[k])$ is zero for all $k$, $1\leq k\leq m$, so $\Hom_{\mathcal D}(C,\bar T[j])$ must be nonzero for some $j$, $1\leq j\leq m$, since $\bar T$ is assumed to be maximal.

By Lemma \ref{nullilokalisert} any such nonzero map must factor through an object in $\mathcal M$ since $\Hom_{\mathcal D_{\mathcal M}}(\widetilde C, \widetilde T[j])=0$. Therefore there must be a nonzero map from $C$ to $M[i]$ as well as a nonzero map from $M[i]$ to $T[j]$ for some $i$. We know that $\Hom_{\mathcal D}(M[i],T[j]) \simeq \Hom_{\mathcal D}(M,T[j-i])$. Since $\Hom_{\mathcal D}(M, T[k])$ is assumed to be zero for $1\leq k \leq m$, the only way for $\Hom_{\mathcal D}(M, T[k])$ to be nonzero for higher $k$ is if $M$ has degree $m$. But then it would be induced by a projective module and so cannot map to higher shifts of $T$ than $m$. This means that $\Hom_{\mathcal D}(M[i],T[j])$ can only be nonzero if $j \leq i$, and so we can assume that $\Hom_{\mathcal D}(C,M[i])$ is nonzero for some $i$ greater than or equal to $1$.

Since $C$ is rigid, we know by Proposition \ref{lagtriangel} that there exists a nonzero indecomposable object $X \in \mathcal D$ such that $\Hom_{\mathcal D}(M,X[t])=0=\Hom_{\mathcal D}(X,M[t])$ for $t=1, 2, \ldots, m$, and $\widetilde{C}\simeq \widetilde{X}$. Therefore $\widetilde T \coprod \widetilde X$ is $m$-rigid in $\mathcal D_{\mathcal M}$, and so any nonzero map from $T$ to $X[r]$ or from $X$ to $T[r]$ for any $r$, $1 \leq r \leq m$, must factor through an object in $\mathcal M$. But Proposition \ref{lagtriangel} says that $\Hom_{\mathcal D}(X,M[t])=0$ for all $t$, so there can be no nonzero maps from $X$ to $T[r]$ for any $r$. We also have $\Hom_{\mathcal D}(M,X[t])=0$ for $t \neq 1-i$, so any nonzero map from $T$ to $X[r]$ must map through some object in $\add \{M[i-1+r]\}$. But $i-1+r \geq 1$ since $i$ and $r$ are assumed to be greater than or equal to $1$, so then we would get a nonzero map from $T$ to a positive shift of $M$. By the same argument as in the previous paragraph we know that $\Hom_{\mathcal D}(T,M[k])=0$ for all $k\geq 0$. We therefore see that $\Hom_{\mathcal D}(\bar T,X[r])=0=\Hom_{\mathcal D}(X,\bar T[r])$ for $r=1, 2, \ldots, m$.

We also get from Proposition \ref{lagtriangel} that $\deg(X)$ must be either $j$ or $j-1$, where $j$ is the degree of $M[i]$ such that $\Hom_{\mathcal D}(C, M[i])\neq 0$. Since $\deg(\widetilde{C})\leq m$ and $C \in \mathcal D_0$, then $\deg(C)\leq m$. This means that $\deg(X) \leq m$ since $\Hom_{\mathcal D}(X,C)\neq 0$. Furthermore $\deg(X)\geq 0$ since $i \geq 1$ and there is a map from $M[i-1]$ to $X$ by Proposition \ref{lagtriangel}. If $\deg(X)=m$, we need to show that $X$ is isomorphic to $P[m]$ for some projective $H$-module $P$. In this case $\deg(C)=\deg(\widetilde{C})$ must also be $m$, and $\widetilde{C}$ is projective in $\mathcal D_{\mathcal M}$.

Consider $\Hom_{\mathcal D}(X,Y[1])$, where $Y$ is an indecomposable object in $(\modu H)[m]$. Since $\widetilde{C}$ is projective, so is $\widetilde{X}$. We also know that $\widetilde{Y[1]}\simeq \widetilde{Y}[1]_{\mathcal D_{\mathcal M}}$, where $[1]_{\mathcal D_{\mathcal M}}$ is the shift functor in $\mathcal D_{\mathcal M}$. Therefore we see that $\Hom_{\mathcal D_{\mathcal M}}(\widetilde{X},\widetilde{Y[1]})=0$ since $\widetilde{X}$ is projective. This means that any nonzero map $X \rightarrow Y[1]$ has to factor through some object in $\mathcal M$. But $\Hom_{\mathcal D}(X,M')=0$ for all $M' \in \mathcal M$ by Proposition \ref{lagtriangel}, and so a map $X \rightarrow Y[1]$ cannot factor through any nonzero $M' \in \mathcal M$. Therefore $\Hom_{\mathcal D}(X,Y[1])$ is zero for any $Y \in (\modu H)[m]$. Hence $X \simeq P[m]$ for some projective $H$-module $P$, and $X \in \mathcal D_G$.

This means that $X \in \add \bar T$ since $\bar T$ is maximal in $\mathcal D_G$. But $\widetilde C \simeq \widetilde X$, and so $X$ being a summand of $\bar T$ means that $\widetilde C$ is a summand of $\widetilde T$. We assumed that $\widehat C \notin \add \widehat T$, so this is a contradiction, and we can conclude that $\widehat T$ is a maximal $m$-rigid object in $\mathcal C_{H'}^m$.

\end{proof}

Before the next theorem we recall the notion of a \emph{complement} of an $m$-rigid object $T'$ in $\mathcal C_H^m$. It is defined as an $m$-rigid object $T''$ which has the properties that no summand of it is in $\add T'$ and that $T' \coprod T''$ is a maximal $m$-rigid object.

\begin{theorem}
Let $T$ be a basic maximal $m$-rigid object in $\mathcal C_H^m$. Then we have the following:

\begin{itemize}
\item [(a)]
$T$ has $n$ indecomposable summands, where $n$ is the number of non-isomorphic simple modules in $\modu H$.
\item [(b)]
If $T'$ is a basic almost complete $m$-rigid object, i.e. has $n-1$ indecomposable nonisomorphic summands, then $T'$ has $m+1$ indecomposable complements.
\end{itemize}

\end{theorem}

\begin{proof}
Both properties obviously hold when $n=1$.

Assume by induction that the claim holds for any finite dimensional hereditary algebra with $n-1$ isoclasses of simple modules, and let $H$ be a finite dimensional hereditary algebra with $n$ isoclasses of simple modules.

First assume that $T$ is a basic maximal $m$-rigid object. If we now localise with respect to any indecomposable summand $M$ of $T$, we get by Theorem \ref{teorem1} that $\mathcal L_{\mathcal M}(T)=\widetilde T$ is a basic maximal $m$-rigid object over some finite dimensional hereditary algebra $H'$ with $n-1$ isoclasses of simples, and so has $n-1$ summands. By Proposition \ref{liktantsummander} we know that nonisomorphic indecomposable summands of $T$ not having $M$ as a summand will correspond to nonisomorphic indecomposable summands of $\widetilde T$. Therefore $T$ must have $n$ summands.

Next assume that $T^*$ is an almost complete basic $m$-rigid object in $\mathcal C_H^m$. By Proposition \ref{liktantsummander} we can localise at any indecomposable summand $M$ of $T^*$, and $\widetilde{T^*}$ will then be an $m$-rigid object in $\modu H' \vee\ (\modu H')[1] \vee \ldots \vee (\modu H')[m-1] \vee H'[m]$ in $D^b(H')$ for some hereditary algebra $H'$ with $n-1$ isoclasses of simple modules. This object will now have $n-2$ nonisomorphic indecomposable summands, and will therefore correspond to an almost complete basic $m$-rigid object in $\mathcal C_{H'}^m$ and so by assumption have $m+1$ complements.

By Theorem \ref{teorem1} any complement $X$ of $T^*$ in $\mathcal D_G$ in $\mathcal D$ will correspond to a nonzero object $\widetilde X$ in $\mathcal D_{\mathcal M}$ that will be a complement of $\widetilde{T^*}$.

We need to check if two nonisomorphic indecomposable complements $X$ and  $X'$ of $T^*$ can correspond to the same object in $\mathcal D_{\mathcal M}$. If they do, there will be some triangle $M' \rightarrow X' \rightarrow X \rightarrow$ in $\mathcal D$ where $M' \in \mathcal M$. This will again mean that there is a nonzero map from an object $M'$ in $\mathcal M$ to $X'$, and from $X$ to $M'[1]$, which means that $M'$ must be zero since $\Hom_{\mathcal D}(M[t],X)=0$ when $t \leq -1$ and $\Hom_{\mathcal D}(X',M[s])=0$ when $s\geq 1$, contradicting $X$ and $X'$ being nonisomorphic.

Next assume that $\widehat Y$ is a complement of $\widehat{T^*}$ in $\mathcal C_{H'}^m$. It corresponds to a complement $\widetilde Y$ of $\widetilde{T^*}$ in $\mathcal D_{\mathcal M}$, and there is a unique object $Y_0$ in $\mathcal D_0$ which corresponds to $\widetilde Y$. If this object is a complement of $T^*$, we are done. If not, $\Hom_{\mathcal D}(Y_0, T^*[i])$ or $\Hom_{\mathcal D}(T^*, Y_0[i])$ is nonzero for some $i$, $1 \leq i \leq m$. But it cannot be the latter since by Lemma \ref{nullilokalisert} any such map would have to factor through an object in $\mathcal M$, and $Y_0 \in \mathcal D_0$. Therefore there must be a nonzero map $Y_0 \rightarrow T^*[i]$ for some $i$, and this map has to factor through some object $M'$ in $\mathcal M$. By assumption $\Hom_{\mathcal D}(M,T^*[k])=0$ for $1 \leq k \leq m$. Also $\Hom_{\mathcal D}(M,T^*[k])=0$ for $k = m+1$, since the degree of the indecomposable summands of $T^*[m+1]$ will be greater than or equal to $m+1$ and if $M$ is of degree $m$ it is projective and so cannot have a nonzero map to $T^*[m+1]$. Obviously $\Hom_{\mathcal D}(M,T^*[k])=0$ for $k > m+1$, and so we see that $\Hom_{\mathcal D}(M,T^*[k])=0$ for all $k \geq 1$. This means that if $\Hom_{\mathcal D}(M[r],T^*[i]) \simeq \Hom_{\mathcal D}(M,T^*[i-r])$ is nonzero for some $i \geq 1$, then $i-r \leq 0$ and so $1 \leq i \leq r$. In other words, $M'$ must have some summand that is a positive shift of $M$ and so there is a nonzero map from $Y_0$ to $M[r]$ for some $r \geq 1$. 

Now by Proposition \ref{lagtriangel} there exists an object $Y_1$ such that $\widetilde Y_1 = \widetilde Y_0$ and in particular $\Hom_{\mathcal D}(M,Y_1[i])= 0 = \Hom_{\mathcal D}(Y_1,M[i])$ for $i = 1, 2, \ldots, m$. We claim that $Y_1$ is a complement to all of $T^*$. Any nonzero map from $Y_1$ to $T^*[i]$ or from $T^*$ to $Y_1[i]$ must factor through some object in $\mathcal M$ by Lemma \ref{nullilokalisert}. By Proposition \ref{lagtriangel}, $\Hom_{\mathcal D}(Y_1,M[t])=0$ for all $t \in \mathbb Z$, and so there cannot be any nonzero maps from $Y_1$ to $T^*[i]$ for any $i$. 

By Proposition \ref{lagtriangel} we also know that $\Hom_{\mathcal D}(M,Y_1[t])=0$ for $t \neq 1-r$. This means that any nonzero map from $T^*$ to $Y_1[i]$ must factor through some object in $\add \{M[r-1+i]\}$, but $r-1+i \geq 1$ since $i$ and $r$ are assumed to be greater than or equal to $1$. This would then give us a nonzero map from $T^*$ to a positive shift of $M$. But $\Hom_{\mathcal D}(T^*,M[k])=0$ for $1 \leq k \leq m$ by definition, and it can only be nonzero for $k>m$ if $T^*$ has summands of degree $m$. But any summand of $T^*$ of degree $m$ will be induced by a projective module and thus cannot map to objects of degree greater than $m$. Therefore $\Hom_{\mathcal D}(Y_1,T^*[i])=0$ for $i=1, 2, \ldots, m$, and $Y_1$ must be a complement of $T^*$. 

So any two nonisomorphic indecomposable complements of $T^*$ will correspond to two nonisomorphic indecomposable complements of $\widetilde{T^*}$, and any indecomposable
complement of $\widetilde T^*$ can be lifted to an indecomposable complement of $T^*$ in $\mathcal D$. Thus $T^*$ and $\widetilde{T^*}$ must have the same number of complements.
\end{proof}

Next we will use induction to show that all maximal $m$-rigid objects in $\mathcal C_H^m$ also are $m$-cluster tilting objects.

\begin{theorem}
$\bar T$ is a maximal $m$-rigid object in $\mathcal C^m_H$ if and only if it is also an $m$-cluster tilting object.
\end{theorem}

\begin{proof}
It is clear that any $m$-cluster tilting object will also be a maximal $m$-rigid object. So assume that $\bar T$ is a maximal $m$-rigid object in $\mathcal C^m_H$, and that $\Ext^i_{\mathcal C^m_H}(T,X)=0$ for $i=1, 2, \ldots ,m$ for some $X \in \mathcal D_G$. Then $\bar T$ corresponds to a maximal $m$-rigid object $T$ in $\mathcal D_G$ by Theorem \ref{erimskift}.

Let $M$ be an indecomposable summand of $T$ such that $\Hom_{\mathcal D}(M,T/M)=0$. Such a summand exists by the following argument: Consider the highest degree of any indecomposable summand of $T$, and let $T'$ be the largest summand of $T$ with this degree. In particular $\Ext^1_{\mathcal D}(T',T')=0$, so $T'$ will be induced by a partial tilting module over $H$. Then $T'$ has a summand $M'$ with the property that $\Hom_{\mathcal D}(M',T'/M')=0$ (see for instance Corollary 4.2 of \cite{HR}), and so also $\Hom_{\mathcal D}(M',T/M')=0$ since the degree of all indecomposable summands of $T$ is less than or equal to the degree of $M'$.

We now assume by induction that $\mathcal L_{\mathcal M}(T)=\widetilde T$ is an $m$-cluster tilting object in the natural fundamental domain of $\mathcal C_{H'}^m$ in $\mathcal L_{\mathcal M}(\mathcal D) = \mathcal D_{\mathcal M}$.

Since $\Ext^i_{\mathcal D}(T, X)=0$ for $i=1, 2,\ldots, m$, also $\Ext^i_{\mathcal D_{\mathcal M}}(\widetilde T, \widetilde X)$ will be zero for each such $i$, and so $\widetilde X$ must be in $\add \widetilde T$.

Now let $T_X$ be in $\add T$ such that $\widetilde X \simeq \widetilde{T_X}$ in $\mathcal D_{\mathcal M}$. Then there is a triangle $X \rightarrow T_X \rightarrow M_X \rightarrow$ in $\mathcal D$ such that $M_X \in \mathcal M$ since $T_X \in \mathcal D_0$. Since $T$ is an $m$-rigid object, $M_X$ cannot have summands of the form $M[i]$ when $i=1, 2, \ldots, m$. Similarly we have that $\Hom_{\mathcal D}(M[-i],X)=0$ for $i=1, 2,\ldots, m$, which excludes $M[-i+1]$ for these $i$ (and in particular $M$) as summands. If $T_X$ has no summands in $H[m]$, then clearly $\Hom_{\mathcal D}(T_X,M[k])=0$ when $k>m$. If some summand of $T_X$ is in $H[m]$, it is induced by a projective module over $H$ and so it cannot map to $M[k]$ when $k>m$. Therefore $M_X$ must be zero and $X \in \add T$.
\end{proof}

Now we move on to look at the \emph{$m$-cluster-tilted algebra $\Gamma=\End_{\mathcal C^m_H}(T)^{\op}$} for an $m$-cluster-tilting object $T$ in $\mathcal C_H^m$. In Section 2 of \cite{BMR} it is shown that if $\Gamma$ is a cluster-tilted algebra for some cluster tilting object in a cluster category, then so is $\Gamma/\Gamma e\Gamma$, where $e$ is the idempotent corresponding to one of the vertices of $\Gamma$. This can be generalized to the $m$-cluster category case.

\begin{theorem}
Let $\Gamma$ be as above, and let $\Gamma e \simeq \Hom_{\mathcal C_H^m}(T,M)$ where $M$ is an indecomposable summand of $T$. Then there is a natural isomorphism $\Gamma/\Gamma e\Gamma\simeq\End_{\mathcal C^m_{H'}}(\widehat T)^{op}$.
\end{theorem}

\begin{proof}
The proof is a direct analogue of the proof of Theorem 2.13 in \cite{BMR}. We will use the notation $\mathcal D' = D^b(H')$, where $H'$ is the hereditary algebra associated with $\mathcal U_X$ where $X$ is the $H$-module such that $M \simeq X[k]$ for some $k$. We know that $\mathcal D'$ is equivalent to $\mathcal D_{\mathcal M}$ (see Section 2 of \cite{BMR} for details).

Letting $\Hom(X,Y)/(Z)$ denote the $\Hom$-space modulo maps factoring through some object in $\add \{Z[t] | t \in \mathbb Z\}$ for some object $Z$, we can write $\Gamma/\Gamma e \Gamma \simeq \Hom_{\mathcal D}(T,T) \coprod \Hom_{\mathcal D}(T,GT)/(M \coprod GM)$. So we need to consider $\Hom_{\mathcal D}(T_a,T_b)/(M \coprod GM)$ and $\Hom_{\mathcal D}(T_a,GT_b)/(M \coprod GM)$ for all indecomposable summands $T_a$ and $T_b$ of $T$.

If $T_X$ is an indecomposable direct summand of $T$ not isomorphic to $M$, we can consider the triangle 
\[M_X \rightarrow T_X \rightarrow U_X \rightarrow\] 
induced by a minimal right $\mathcal M$-approximation of $T_X$ in $\mathcal D$. By Lemma \ref{blirnull} we know that $U_X \in \mathcal D_0$, and so by Lemma \ref{dnullogtau} we get $\widetilde{G_{\mathcal D}U_X} \simeq \widetilde{\tau_{\mathcal D}^{-1}U_X[m]_{\mathcal D}} \simeq \tau_{\mathcal D'}^{-1} \widetilde{U_X[m]_{\mathcal D}}\simeq \tau^{-1}_{\mathcal D'}\widetilde{U_X}[m]_{\mathcal D'}\simeq G_{\mathcal D'}\widetilde{U_X}$, where $[1]_{\mathcal D}$ denotes the shift functor in $\mathcal D$ and $[1]_{\mathcal D'}$ the shift functor in $\mathcal D'$.

By Theorem \ref{erimskift} we can assume that $T \in \modu H \vee (\modu H)[1] \vee \ldots \vee (\modu H)[m-1]$. Now pick two indecomposable direct summands $T_a$ and $T_b$ of $T$, and consider the triangles 
\[M_b \rightarrow T_b \rightarrow U_b\rightarrow\] 
and 
\[GM_b \rightarrow GT_b \rightarrow GU_b\rightarrow\] 
induced by the minimal right $\mathcal M$-approximation $M_b\rightarrow T_b$. Then $U_b$ will be in $\mathcal D_0$ by Lemma \ref{blirnull}. We also note that since $M_b \rightarrow T_b$ is a minimal right $\mathcal M$-approximation, $GM_b \rightarrow GT_b$ will be a right $G\mathcal M$-approximation, where by $G\mathcal M$ we mean the full subcategory of $\mathcal D$ with objects $\{GX | X \in \mathcal M\}$.

If we apply $\Hom_{\mathcal D}(T_a,-)$ to the last triangle, we get the long exact sequence 
\begin{align*}
	\Hom_{\mathcal D}(T_a, GM_b)	& \rightarrow \Hom_{\mathcal D}(T_a,GT_b) \\
																&	\rightarrow \Hom_{\mathcal D}(T_a,GU_b)	\rightarrow	\Hom_{\mathcal D}(T_a,GM_b[1]).
\end{align*}
The last term of the sequence is zero since the degree of $GM_b[1]$ is greater than or equal to $m+1$, and if $T_a$ has degree $m$ it is induced by a projective $H$-module and thus any map to an object in the $(m+1)$-th shift is zero. We therefore see that $\Hom_{\mathcal D}(T_a,GU_b) \simeq \Hom_{\mathcal D}(T_a,GT_b)/(GM)$.

Next we claim that there is an exact sequence
\begin{align*}
	\Hom_{\mathcal D}(T_a,GM_b)/(M)	& \rightarrow \Hom_{\mathcal D}(T_a,GT_b)/(M) \\
																	&	\rightarrow \Hom_{\mathcal D}(T_a,GU_b)/(M) \rightarrow 0
\end{align*}
induced by the above long exact sequence. For the argument for this we refer the reader to the analogous argument in \cite{BMR} (p. 156-157). It follows from this that $\Hom_{\mathcal D}(T_a,GU_b)/(M) \simeq \Hom_{\mathcal D}(T_a,GT_b)/(M\coprod GM)$. 

Now we form the triangle 
\[M_1 \rightarrow GU_b \rightarrow (GU_b)'\rightarrow\]
where $M_1 \rightarrow GU_b$ is the minimal right $\mathcal M$-approximation of $GU_b$. Again we apply $\Hom_{\mathcal D}(T_a,-)$ to get the long exact sequence
\begin{align*}
	\Hom_{\mathcal D}(T_a,M_1)	& \rightarrow \Hom_{\mathcal D}(T_a, GU_b) \\
															&	\rightarrow \Hom_{\mathcal D}(T_a, (GU_b)') \rightarrow \Hom_{\mathcal D}(T_a,M_1[1]).
\end{align*}
We see that $\Hom_{\mathcal D}(T_a,M_1[1])=0$ by the following argument: 

Since there is a nonzero map from $T_b$ to $U_b$ induced by a minimal right $\mathcal M$-approximation, we have $0 \leq \deg(U_b)$. This means that $m \leq \deg(GU_b)$. Since $M_1 \rightarrow GU_b$ is a minimal right $\mathcal M$-approximation, all nonzero summands of $M_1$ have a nonzero map to $GU_b$. Therefore the degree of the indecomposable summands of $M_1$ is greater than or equal to $m-1$, and so the degree of the indecomposable summands of $M_1[1]$ is greater than or equal to $m$. Since $M$ is assumed to be of degree less than or equal to $m-1$, the indecomposable summands of $M_1[1]$ must be positive shifts of $M$. This means that $\Hom_{\mathcal D}(T_a,M_1[1])=0$ since $T$ is $m$-rigid. If the map $M_1 \rightarrow GU_b$ is zero, $M_1$ must be zero since it is induced by the minimal right $\add M$-approximation of $GU_b$, and so obviously $\Hom_{\mathcal D}(T_a,M_1[1])=0$.

We know that a map $T_a \rightarrow GU_b$ will factor through an object in $\mathcal M$ if and only if it factors through the minimal right $\mathcal M$-approximation of $GU_b$. This means that by the previous long exact sequence we get that $\Hom_{\mathcal D}(T_a,GU_b)/(M) \simeq \Hom_{\mathcal D}(T_a,(GU_b)')$. Since $(G_{\mathcal D}U_b)'$ and $U_b$ are in $\mathcal D_0$, we have 
\begin{align*}
	\Hom_{\mathcal D}(T_a, (G_{\mathcal D}U_b)')	& \simeq \Hom_{\mathcal D'}(\widetilde{T_a}, \widetilde{(G_{\mathcal D}U_b)'}) \simeq \Hom_{\mathcal D'}(\widetilde{T_a}, \widetilde{G_{\mathcal D}U_b}) \\
																								& \simeq \Hom_{\mathcal D'}(\widetilde{T_a},G_{\mathcal D'}\widetilde{U_b}) \simeq \Hom_{\mathcal D'}(\widetilde{T_a},G_{\mathcal D'}\widetilde{T_b}) &
\end{align*}
by Proposition \ref{dnulliso}. This means that there is an isomorphism $\Hom_{\mathcal D}(T_a, G_{\mathcal D}T_b)/(M \coprod GM) \simeq \Hom_{\mathcal D}(T_a,GU_b)/(M) \simeq \Hom_{\mathcal D'}(\widetilde{T_a}, G_{\mathcal D'}\widetilde{T_b})$. 

We have that $\Hom_{\mathcal D}(T_a,T_b)/(M \coprod GM)$ is isomorphic to $\Hom_{\mathcal D}(T_a,T_b)/(M)$ since $T_a$, $T_b$ and $M$ are in $\mathcal D_G$ and so no map from $T_a$ to $T_b$ can factor through $GM$. By applying $\Hom_{\mathcal D}(T_a,-)$ to the triangle $M_b \xrightarrow{f_b} T_b \rightarrow U_b$ in $\mathcal D$ induced by the minimal right $\mathcal M$-approximation of $T_b$, we get the exact sequence
\[\Hom_{\mathcal D}(T_a,M_b) \rightarrow \Hom_{\mathcal D}(T_a,T_b) \rightarrow \Hom_{\mathcal D}(T_a,U_b) \rightarrow \Hom_{\mathcal D}(T_a,M_b[1]).\]
We have that the degree of the indecomposable summands of $M_b$ is less than or equal to $m-1$ since $M_b \rightarrow T_b$ is a minimal right $\mathcal M$-approximation. This implies that $\Hom_{\mathcal D}(T_a,M_b[1]) = 0$. Therefore $\Hom_{\mathcal D}(T_a,T_b)/(M) \simeq \Hom_{\mathcal D}(T_a,U_b) \simeq \Hom_{\mathcal D'}(\widetilde{T_a},\widetilde{T_b})$, the latter isomorphism coming from Proposition \ref{dnulliso}. This again means that $\Hom_{\mathcal D}(T_a,T_b)/(M \coprod GM)\simeq \Hom_{\mathcal D'}(\widetilde{T_a},\widetilde{T_b})$. Putting everything together we see that there is a vector space isomorphism $\Gamma / \Gamma e \Gamma = \Hom_{\mathcal D}(T,T) \coprod \Hom_{\mathcal D}(T,GT)/(M \coprod GM) \simeq \Hom_{\mathcal D'}(\widetilde T,\widetilde T) \coprod \Hom_{\mathcal D'}(\widetilde T,G_{\mathcal D'}\widetilde T)$, as desired. One can check that this isomorphism also is a ring isomorphism, and so we are done. 

\end{proof}

\begin{corollary}
Let $\Gamma$ be an $m$-cluster tilted algebra and let $e$ be the idempotent corresponding to a vertex of $\Gamma$. Then $\Gamma/\Gamma e\Gamma$ is also $m$-cluster tilting.
\end{corollary}

%%%%%%%%%%%%%%%%%%%%%%%%%%%%%%%%%%%%%%%%%%%%%%%%%%%%%%%%%%%%%%%%%%%%%%%%%%%%%%%

\end{document}